\documentclass[12pt]{article}

\title{The SIR model of an epidemic}
\author{William G. Faris\\Department of Mathematics\\University of Arizona}
\date{August 9, 2020}

\usepackage{graphicx}

\newtheorem{proposition}{Proposition}

\newcommand{\sech}{\mathrm{sech}}

\begin{document}
\maketitle

\begin{abstract}

The SIR model is a three-compartment model of the time development of an epidemic. After normalizing the dependent variables, the model is a system of two non-linear differential equations for the susceptible proportion $S$ and the infected proportion $I$.  After normalizing the time variable there is only one remaining parameter. This largely expository article is mainly about aspects of this model that can be understood with calculus. It also discusses an alternative exactly solvable model that appeared in early work of Kermack and McKendrick.  This model may be obtained by replacing $SI$  factors  by $\sqrt{2S-1}I$ factors. For a mild epidemic,  where $S$ is decreasing from 1 but remains fairly close to 1,  this is a reasonable approximation. 
\end{abstract}

\section{The simplest model of an epidemic}

The  SIR model  is the simplest differential equation model that describes how an epidemic begins and ends. It depends on only two parameters: One  governs the timing,   the other 
 determines everything else. It gives a glimpse into the world of more complicated epidemic models.  
 
 The SIR model is standard in the literature of epidemiology \cite{H,KM}, and it even shows up in textbooks on calculus and differential equation. There is no exact formula for how the infection level depends on time. This account reviews what can be done in spite of this limitation. 
 There is a  modification of the SIR model for which there is an exact formula for the time dependence. 
 This account also presents a way of understanding this modified model .

 For a fixed  population $S$ is the proportion of susceptibles,  $I$ is the proportion of infected, and $R$ is the proportion of removed. Removal is supposed to be permanent. The mechanism of removal is not specified. Some individuals might not survive; others  might recover with full immunity.

Every  member of the population is supposed to belong  to one of these three compartments,  so $S + I + R = 1$.
While S, I, and R change with time, the parameters that define the model do not change with time. This is a closed system, with no reintroduction from the outside. 

There are two parameters in the model. One is a rate $b$ that governs how fast the susceptible-infected pairs are changing into infected-infected pairs. The other is a rate $a$ that governs how fast infected individuals are being removed. 
The ratio 
\begin{equation}
r_0 = \frac{b}{a} 
\end{equation}
is the \emph{basic reproduction number} that determines the scope of the epidemic.  This represents the  number of new infections per removal at the beginning of the epidemic. 
When $r_0<1$ the epidemic immediately begins to die out. The more interesting situation is when $r_0 > 1$. 

In this model there are no interventions; the rate constants $b$ and $a$ never change, and so also the reproductive number $r_0$ is constant. 
The goal is to see how the course of the epidemic depends on the parameters. 

Here is an outline of what happens. During the course of the epidemic $S$ decreases and $R$ increases. What happens to $I$ is governed by the \emph{effective reproduction number} $r_0 S$. This represents the average number of new infections per removal. As the pool of susceptibles is depleted, 
this number decreases with time. As long as $r_0 S$ is greater than one, the proportion $I$ of infected increases. It reaches a peak when $r_0 S = 1$. Afterward  $r_0 S$ is less than one, and $I$ decreases to zero. 
Meanwhile, $R$ increases to a final value representing the proportion of the population  eventually  infected and removed.   
 
 It is sometimes useful to take the fundamental 
parameters as $r_0$ and $a$. In that case, the role of $r_0$ is to describe the size of the epidemic, 
while the role of $a$ is to describe how fast the epidemic takes place as a function of time. 

Figures 1, 2, and 3 show the behavior of the SIR model for $r_0$ values 2, 3, and 6. The following discussion will clarify what these pictures mean. Even when $r_0$ is 2, the final value of R is 80 percent.   
Also, the $r_0$ equal to 6 picture shows that the infection proportion curve $I$ is not symmetric about its peak. 
This is because at this level of $r_0$  infection is fast and removal is slow. (The modified SIR model will show why
for small $r_0$ the infection curve is approximately symmetric.)

Figure~4 shows the dependence of the final proportion $R$ of removed on the basic reproduction number $r_0>1$. The values are shockingly large; the growth of infection has a dramatic effect before it is checked.

\section{The SIR model}

\begin{figure}[t]

\includegraphics*{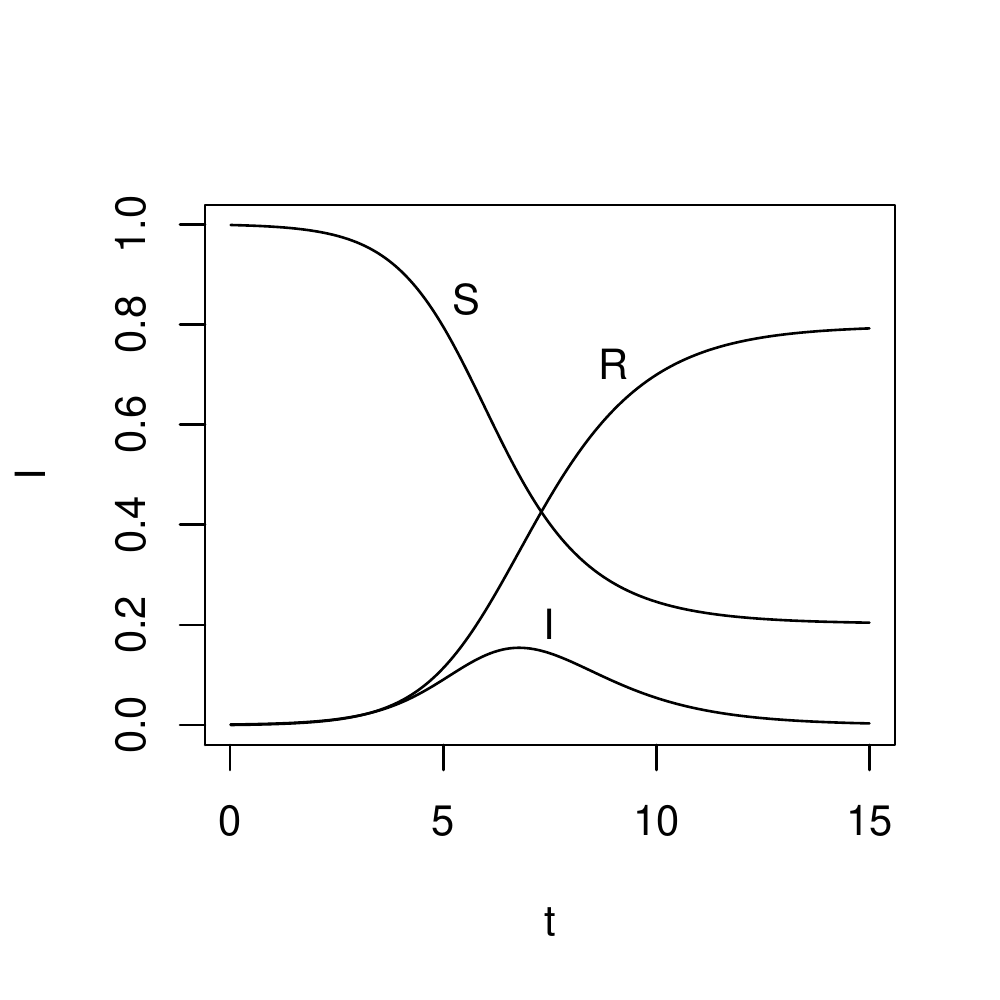}
\caption{$S,I,R$ as functions of $\tau = at$: $r_0 = 2$.}
\end{figure}

The SIR model has fundamental equations
\begin{eqnarray}
\frac{dS}{dt}  &=& - b  SI,  \nonumber \\
\frac{dI}{dt}  &=&  b SI - a I \nonumber \\
\frac{dR}{dt} &=& a I .
\end{eqnarray}
Here $S>0$ is the proportion susceptible, $I> 0$ is the proportion   infected, $R\geq0$ is the proportion removed. 
Since $S + I +R = 1$, any two of these equations implies the third one. 
There are two rate constants: $b>0$ and $a > 0$. The initial conditions at $t=0$ are $S_0 <1 $ and $I_0 >  0$ with $S_0 + I_0 = 1$. This implies that $R_0 = 0$. 
It is often useful to think of $S_0 \approx 1$ and $I_0 \approx 0$ as being close to 1 and 0.

It is convenient to define  a rescaled time $\tau = a t$. This is measuring time in units that roughly correspond to the average time for removal. Since the basic reproduction number $r_0 = b / a$, the equations become
 \begin{eqnarray}
\frac{dS}{d\tau}  &=& -  r_0 SI, \\
\frac{dI}{d\tau}  &=&  (r_0 S - 1) I, \\
\frac{dR}{d\tau} &=& I .
\end{eqnarray}
 The proportion  of susceptibles $S$ starts at $S_0$  and steadily decreases, eventually approaching a limit $S_\infty>0$ as $\tau \to+\infty$. 
 The proportion of  removeds  $R$ starts at 0 and steadily increases, eventually approaching a limit $R_\infty<1$ as $\tau \to +\infty$. At each time the growth rate of the proportion $I$ of infecteds is approximately exponential with rate $r_0 S - 1$. 
 In the long run this is negative, and 
 the proportion of infecteds $I$ approaches 0 as $\tau \to +\infty$.  It follows that $S_\infty + R_\infty = 1$.

By definition the effective reproduction number is $r_0 S$. From the equation 
\begin{equation}
-\frac{dS}{dR} = r_0 S 
\end{equation}
it  appears that the effective reproduction number represents new infections per removal. 

Since $I = 1 - S - R$, it follows that
\begin{equation}
\frac{dI}{dR}  = r_0 S - 1.
\end{equation}
This measures net infections  per removal.   
This equation for $dI/dR$ leads to the first and most basic result. Since $dI/d\tau = dI/dR \; dR/d\tau$ and $dR/d\tau > 0$, it can
be thought of as a result either about dependence on $R$ or dependence on $\tau$. 
Suppose $r_0 S_0 > 1$. 
\begin{itemize}
\item The infection proportion $I$ as a function of the removed proportion $R$ starts at $I_0$, increases to a peak at the $R$ level where the effective reproduction number $r_0 S = 1$, and then decreases to zero at $R_\infty$.
\item The infection proportion $I$ as a function of time $\tau$ starts at $I_0$, increases to  a peak at the time $\tau$ when the effective reproduction number $r_0 S = 1$, and then decreases to zero as $\tau \to +\infty$.
\end{itemize}

\begin{figure}[t]

\includegraphics*{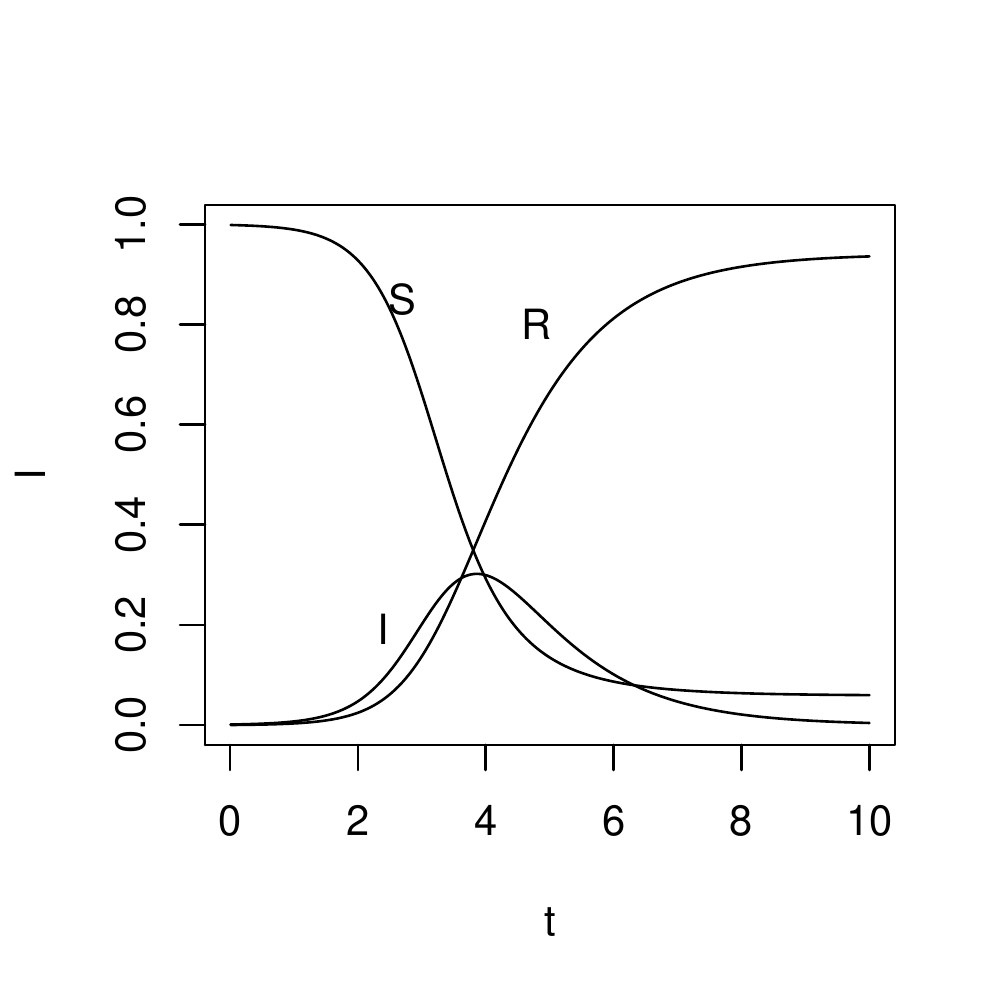}
\caption{$S,I,R$ as functions of $\tau=at$: $r_0 = 3$.}

\end{figure}
 
\section{$S$ as a decreasing function of $R$}

The differential equation $dS/dR = - r_0 S$ is the equation for exponential decay. 
If  $R_0 = 0$ this has solution 
\begin{equation}
S = S_0 e^{- r_0 R} 
\end{equation}
for $0 \leq R \leq R_\infty$. 

The $S$ and $R$ values determine each other. In fact,
\begin{equation}
R = - \frac{1}{r_0} \ln\left( \frac{S}{S_0} \right) 
\end{equation}
for $S_0 \geq S \geq S_\infty$. 

The infection is then 
\begin{equation}
I  = 1 - S - R = 1 - S_0 e^{-r_0 R}  - R. 
\end{equation}
Alternatively, 
\begin{equation}
I =  1 - S - R = 1 - S  +  \frac{1}{r_0}  \ln \left( \frac{S}{S_0} \right).
\end{equation}
Its values  are positive; they start at $I_0 = 1 - S_0$ and end at zero as $S$ arrives at $S_\infty$. 

Putting these explicit formulas back in the original equations leads to integral expressions for $t$ as a function of $R$ and for $t$ as a function of $S$.
See for instance \cite{HLM}.

\begin{figure}[t]

\includegraphics*{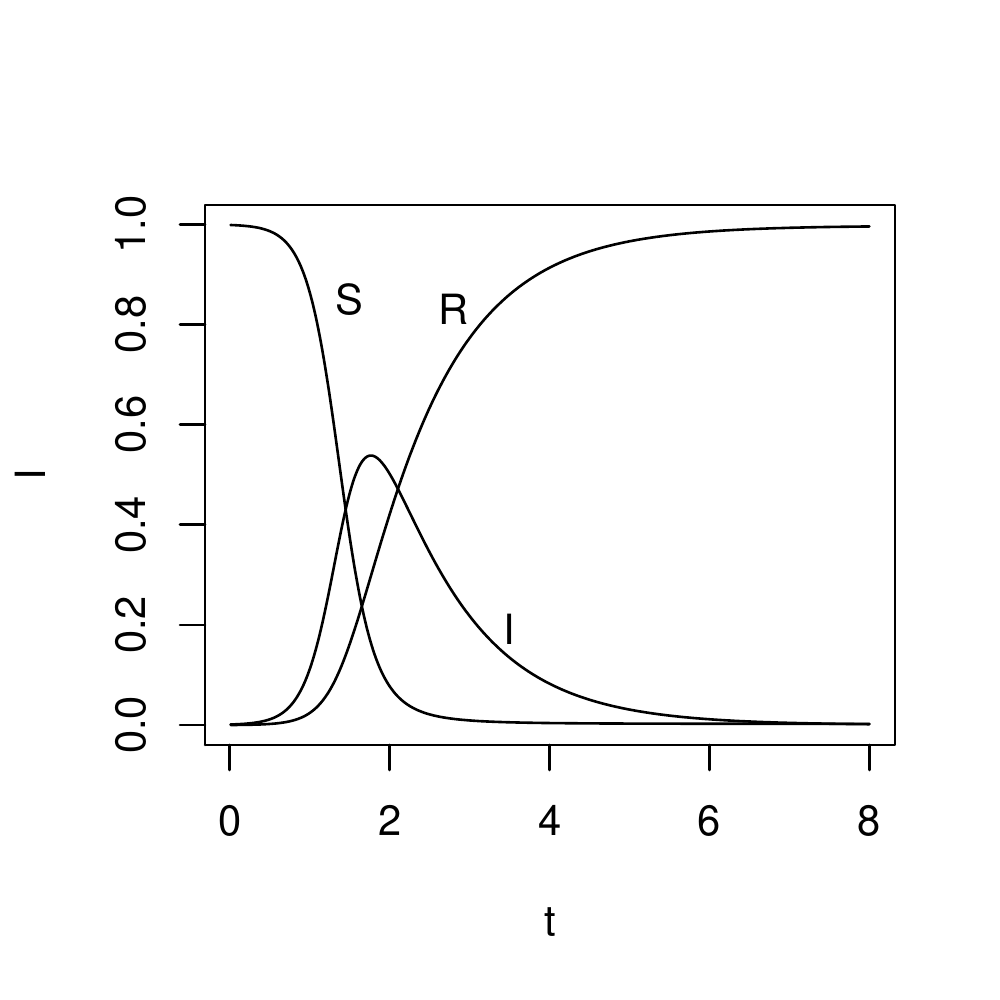}
\caption{$S,I,R$ as functions of $\tau = at$: $r_0 =6$}

\end{figure}

\section{$I$ at peak infection}

\begin{proposition}
Start with $r_0 S_0 > 1 $. The peak infection takes place when the effective reproduction number $r_0 S = 1$.
At peak infection the values of $S, I , R$ are 
\begin{eqnarray}
S^* &=& \frac{1}{r_0} , \nonumber \\
I^*  &=& 1   - \frac{1}{r_0  }  + \frac{1}{r_0 } \ln \left( \frac{1}{r_0 S_0} \right) \nonumber \\
R^* &=& - \frac{1}{r_0 } \ln \left( \frac{1}{r_0 S_0} \right)  .
\end{eqnarray}
As $r_0$ increases from $1/S_0$ to $+\infty$,  the   $S^*$ value  decreases from $S_0$ to 0, the $I^*$  value increases from $I_0$ to  1,
and the   $R^*$ value rises from 0 to $S_0/e$ (when $r_0 = e/S_0)$ and then falls back to 0. 
\end{proposition}


With the approximation $I_0 \approx 0$ and $S_0 \approx 1$ these peak $I$ values depend only on $r_0$. 
To calculate how big the epidemic gets, all one needs is the value of $r_0$ and a scientific calculator. 

At peak $I$ the value of $R$ can never exceed $1/e \approx 0.3679$. The remaining growth of $R$ is after the peak. When $r_0$ is very large, then $R_\infty$
is close to one, and almost all the growth of $R$ from 0 to $R_\infty$ is after the peak. The onset is sudden; the effects linger.

\section{$S$ and $R$ at the end of the epidemic}

\begin{proposition}
During the epidemic  $S$ decreases to $S_\infty$, $I$ eventually decreases to 0, and $R$ increases to $R_\infty$. The limiting values $S_\infty$ and  $R_\infty$ 
satisfy $S_\infty + R_\infty = 1 $, and
\begin{equation}
S_\infty = S_0 e^{-r_0 R_\infty}.
\end{equation}
Furthermore, the equation 
\begin{equation}
1- R_\infty = S_0 e^{-r_0 R_\infty}
\end{equation}
 has a unique solution with $0<R_\infty < 1$.  Increasing $r_0$ makes 
 this solution closer to one. 
\end{proposition}

When $S_0 \approx 1$ this equation takes the approximate form
\begin{equation}
1- R_\infty \approx  e^{r_0 R_\infty}
\end{equation}
There is no explicit formula for the solution. However at peak $I$ the value of $S$ is $1/r_0$, so $S_\infty < 1/r_0$. It follows that $1 - 1/r_0 < R_\infty$. 

The equation is  a fixed point equation for $R_\infty$. It may be solved by iteration. The iteration function is $g(R) = 1 - e^{-r_0 R}$.  Start with the value $R =  1 - 1/r_0$. 
Keep replacing $R$ by $g(R)$.  These iterations give larger and larger 
values of R, increasing to the fixed point $R_\infty$. To calculate the end result of the epidemic, all that one needs is the value of $r_0$, a scientific calculator, and  patience. 

Table~1 shows the values $S^*, I^*, R^*$ at peak infection for various values of $r_0$. It also shows the end values $S_\infty, R_\infty$. For $r_0 = 6$ the proportion $R^*$ removed 
at peak infection is not very large; the onset of the infection is rapid, and the removals have not yet had time to occur.

\begin{table}[t]
\begin{center}
\begin{tabular}{l | l l l | l l}
&$S^*$  & $I^*$ & $R^*$ &  $S_\infty$ & $R_\infty $  \\ \hline
$r_0 = 2$ &0.50 & 0.15 & 0.35 & 0.20 & 0.80 \\
$r_0 = 3$ & $0.33$ & $0.30 $ & $0.37$ & $0.06$ & $0.94$  \\
$r_0 = 6$ &  $0.17$ & $0.53 $ & $0.30$ & $0.003$ & $0.997$ \\
\end{tabular}
\end{center}
\caption{Peak and end values in SIR model}

\end{table}

Figure 4 shows the dependence of $R_\infty$ on $r_0$ for for all values $0 < r_0 \leq 6$. 

\begin{figure}[t]

\includegraphics*{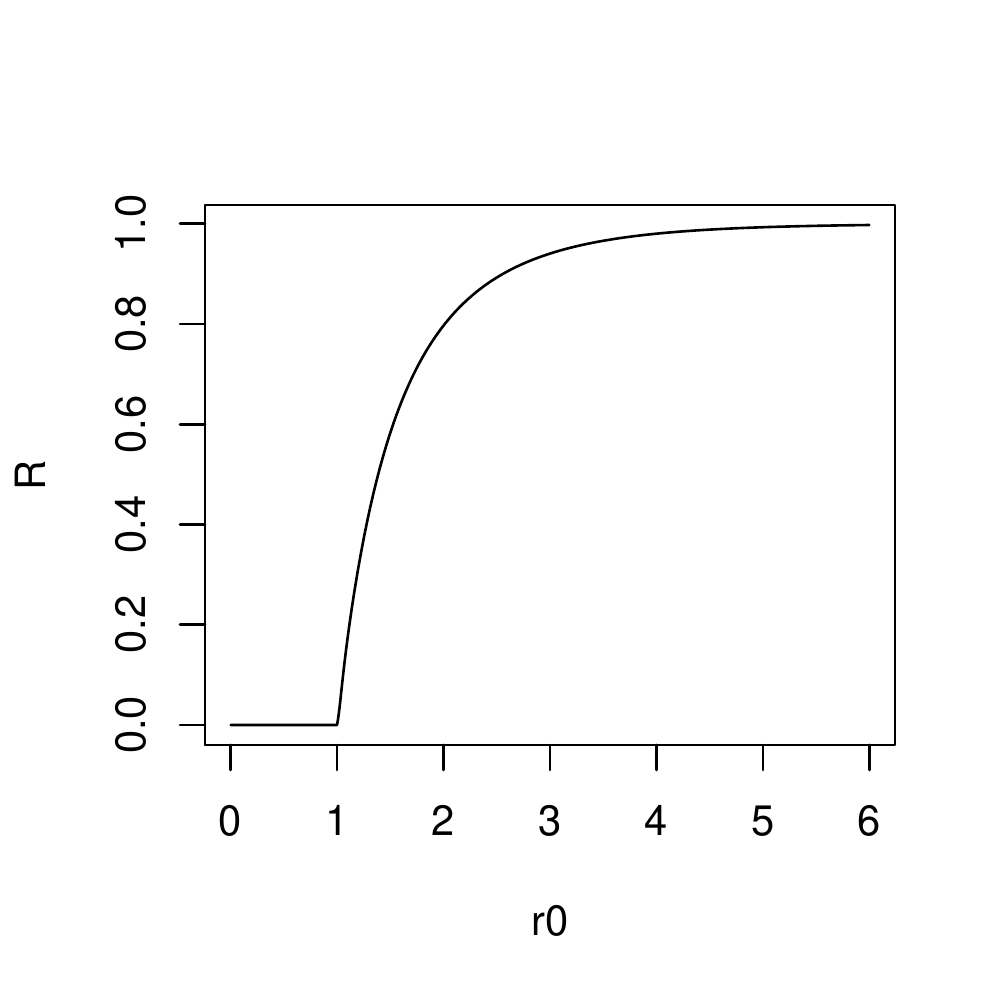}

\caption{$R_\infty$ as a function of $r_0$ for  $r_0 > 1$.}

\end{figure}

 \section{The fastest increase of new infections}
 
 How fast can the proportion of infections increase? There are two obvious measures. The rate of increase of new infections is $-dS/d\tau$. The rate of increase or decrease of infections is $dI/d\tau$. 
 They are related by $dI/d\tau \leq - dS/d\tau$. The results for $-dS/d\tau$ are considerably simpler. 
 
 The most obvious result is an upper bound resulting from $S+I \leq 1$:
  \begin{equation}
-\frac{dS}{d\tau}  = r_0 SI \leq   \frac{1}{4} r_0.
 \end{equation}
Since $\tau = a t$ and $b  = r_0 a$, it follows that
 \begin{equation}
 -\frac{dS}{dt} \leq \frac{1}{4} r_0 a = \frac{1}{4} b.
 \end{equation}
 The constants $b$ and $a$ set the time scale in  usual time units, days or weeks or months.

 A more detailed analysis involves second derivatives.  For $S$ the result is
 \begin{equation}
 \frac{d^2 S}{d\tau^2} =    r_0 SI (r_0 I -   ( r_0 S- 1)  ) . 
 \end{equation}
  Suppose $I_0 \approx 0$ and $S_0 \approx 1$. When $r_0 > 1$, the second derivative of $S$ starts out negative and eventually becomes positive. Thus $S$ has an inflection point where $dS/d\tau$ is minimal.
 
 \begin{proposition}
 The fastest increase of new infections is where $-dS/d\tau$ is largest. 
 This is where 
 \begin{equation}
r_0  I = r_0 S - 1.
 \end{equation}
  The  value is 
 \begin{equation}
-\frac{dS}{d\tau} = r_0 SI =   S (r_0 S - 1). 
 \end{equation} 
 \end{proposition}
 
 The relevant equation 
 \begin{equation}
 r_0 \left(   1 - S  +  \frac{1}{r_0}  \ln \left( \frac{S}{S_0} \right) \right) = r_0 S - 1
 \end{equation}
 describes the intersection of a curve with a line; it has a unique solution. 

 A similar calculation gives equations for the two points where $dI/d\tau$ has maximum or minimum values.

\section{The SI model}

The original $SIR$ model makes sense when $b> 0$ and $a = 0$. This describes an infection that runs through the entire population. Since $R$ stays at zero,
it could be called the $SI$ model. More commonly, it is called the logistic model. It has an explicit solution as a function of $t$. This explicit solution will
be used in the following section for quite another purpose.

In the SI model the variables $S, I$ satisfy $S + I = 1$. The equations  are
\begin{eqnarray}
\frac{dS}{dt}  &=& - b SI,  \nonumber \\
\frac{dI}{dt}  &=&  b SI . \\
\end{eqnarray}
They may be solved explicitly to give $S$ and $I$ as functions of time. 
Since 
\begin{equation}
\frac{1}{b} \frac{d}{dt} \ln\left( \frac{I}{S} \right) = S + I = 1,
\end{equation}
integration gives
\begin{equation}
t =t^*  + \frac{1}{b} \ln \left( \frac{I}{S} \right)  =  t^* + \frac{1}{b} \ln \left( \frac{I}{1-I} \right).
\end{equation}
Here $t^* $ is a constant of integration.
The inverse function is given by
\begin{equation}
I = \frac{1}{2} + \frac{1}{2} \tanh\left( \frac{b (t-t^* )}{2} \right) .
\end{equation}
This is the general solution. It is the cumulative probability function for a logistic probability distribution. The value $t = t^* $ corresponds to $I = \frac{1}{2}$.
In this simple model the infection   grows to take over the entire population. 

The derivative is 
\begin{equation}
\frac{dI}{dt} = \frac{b}{4} \sech^2 \left( \frac{b (t-t^*)}{2} \right).
\end{equation}
This is the density for a logistic probability distribution with mean $t^*$ and standard deviation $\pi/\sqrt{3}$ times $1/b$. 

The relation between the SIR model and the SI model is subtle. As $r_0 \to +\infty$ and $a \to 0$ with $b = r_0 a$ fixed, the  convergence is very non-uniform.  
For each $a>0$ the infection proportion $I$ goes to zero as $t \to +\infty$, while for $a = 0$ the proportion $I$ goes to 1. 
The SI model only captures the onset of the SIR epidemic, well before the peak.

\section{The modified SIR model}

In the SIR model the function $I$ as a function of $t$ is not symmetric about the peak, and there is no obvious formula for this function. 
The paper by  Kermack and McKendrick \cite{KM} introduced an approximation for this function near the peak that is not only symmetric, but
admits an explicit formula for $I$ as a function of $t$ in terms of the square of a hyperbolic secant function. This section presents a modified 
SIR model whose exact solution is precisely the approximation given by Kermack and McKendrick. 
(The paper \cite{M} describes another variant SIR model that also has an explicit formula for the time dependence.)

In the modified model the effective reproduction number is given for $1 \geq S \geq \frac{1}{2}$ by $r_0 \sqrt{1 + 2(S-1)}$,  which is less than $r_0 S$. As the proportion of of susceptibles decreases
from 1, their effective participation in feeding the epidemic also weakens. 
To make the formulas simple we again make the approximation that $S_0 \approx 1$ at time zero.

Again $S + I + R = 1$. Take $r_0 > 1$. The modified SIR equations are  
 \begin{eqnarray}
\frac{dS}{d\tau}  &=& -  r_0 \sqrt{2S-1} \, I, \\
\frac{dI}{d\tau}  &=&  (r_0 \sqrt{2S-1} - 1) I, \\
\frac{dR}{d\tau} &=& I .
\end{eqnarray}
 These equations use the  modified effective reproduction number 
 \begin{equation} 
 r_0 \sqrt{2S - 1} = r_0 \sqrt{ 1 + 2(S-1)} =r_0 S \sqrt{1 - \left( \frac{S-1}{ S } \right)^2 } \leq r_0 S 
 \end{equation}
 for $S \geq 1/2$. For $S$ near 1 this is approximately equal to the effective reproduction number $r_0 S$ of the SIR model. Even for $S = 5/8$  the value of $\sqrt{2S-1}$ is  only down to  $1/2$, which is not so far off. 
   This effective reproduction number could  apply in a  scenario where susceptibles flee as they see their numbers diminish. Its nice properties will appear in 
   the discussion that follows.  
 
It follows that 
\begin{equation}
-\frac{dS}{dR} = r_0 \sqrt{2S-1}.
\end{equation}
and 
\begin{equation}
\frac{dI}{dR}  = r_0 \sqrt{2S-1}  - 1.
\end{equation}
Here the  value of $S$ where $dI/dR = 0$ is $(1/2) ( 1 + 1/r_0^2 )$,  always greater than $1/2$. 

The   differential equation  for $S$ and $R$ has solution
\begin{equation}
r_0 R = 1- \sqrt{2S-1}.
\end{equation}
This can be inverted to express $S$ as a quadratic function. 
\begin{equation}
S = 1 - r_0 R  + \frac{1}{2} r_0^2 R^2. 
\end{equation}
 This is the second order Taylor expansion of the exponential function from the usual SIR model.

\begin{figure}[t]

\includegraphics*{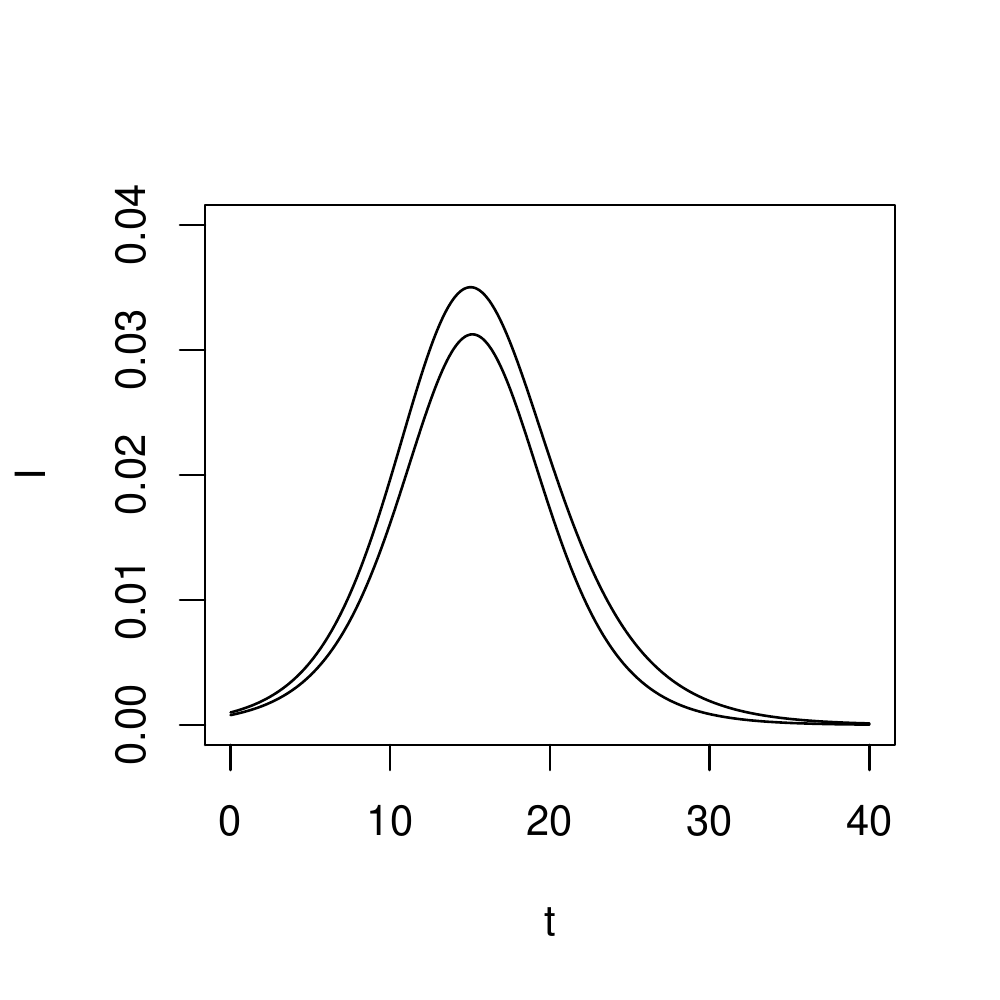}

\caption{$I$ as function of $\tau = at$ for SIR and modified SIR:  $r_0 =4/3$.}

\end{figure}

The time dependence may be found using
 \begin{equation}
\frac{dR}{d\tau} = I =   (r_0 - 1) R - \frac{1}{2} r_0^2 R^2.
\end{equation}
Let
\begin{equation}
R = 2\frac{1}{r_0} \left( 1 - \frac{1}{r_0} \right) \bar R. 
\end{equation}
With this change of variable 
the equation becomes
\begin{equation}
\frac{d\bar R}{d\tau} = (r_0-1) \bar R (1 - \bar R).
\end{equation}
This is the logistic equation with known solution. Define $\tau^*$ as the time when $\bar R = \frac{1}{2}$. 
The solution for $R$ is
\begin{equation}
R = \frac{1}{r_0} \left( 1 - \frac{1}{r_0} \right) \left( 1 + \tanh\left( \frac{1}{2} (r_0-1) (\tau - \tau^* ) \right) \right)
\end{equation}
As $t \to +\infty$ the proportion $R$ increases to the limiting value $2 (1/r_0)(1 - 1/r_0)$. This never exceeds $1/2 $. 
Also
\begin{equation}
I = \frac{dR}{d\tau} = \frac{1}{2}   \left( 1 - \frac{1}{r_0} \right)^2  \sech^2 \left( \frac{1}{2} (r_0-1) (\tau - \tau^* ) \right). 
\end{equation}
The infection profile is given by a the square of a hyperbolic secant function. The maximum is at $\tau^*$. Infection never touches more than half the population. 

This model should not be taken too seriously for large $r_0$. In this case the main part of the infection takes place over a short period of time, roughly $1/(r_0 -1)$,  and it 
has little consequence in the long run. For $r_0$ only a little larger than one the results should be a reasonable approximation to those of the SIR model. At the peak 
the value of the expansion parameter $r_0 R$ is $1 - 1/r_0$ which is close to zero. 

Figure~5 for $r_0 = 4/3$ shows $I$ as a function of $\tau = at$ for the SIR model and for the modified SIR model. They look roughly the same, though the modified SIR model predicts a smaller epidemic. The fit could be made better by 
adjusting the parameters. However the lack of symmetry about the peak in the SIR model will never be reflected in the modified model.





\section{Conclusion}

In the SIR model most essential features are captured in one number, the basic reproduction number. The formulas for peak infection and for long-time behavior depend only on this number. It determines everything about the solution, except for the timing. 

This model exhibits the  basic features of an epidemic. The infection initially rises, but eventually removal is dominant, and the infection dies out. The mechanism is that the proportion of susceptibles decreases until there are too few to maintain the infection. 
If the basic reproduction number is moderately  larger than one,  then the final proportion removed is large (Figure~4). Whether this is benign or sinister depends on the consequences of  removal. 

The author thanks Robert Indik for valuable comments.


\begin{thebibliography}{9}

\bibitem{HLM} Tiberiu Harko, Francisco S. N. Lobo, and  M. K. Mak, 
Exact analytical solutions of the Susceptible-Infected-Recovered (SIR) epidemic model and of the SIR model with equal death and birth rates,
Applied Mathematics and Computation 236 (2014),  184--194.

\bibitem{H}  Herbert W. Hethcote, Three basic epidemiological models,  pp 119-144,     in  Applied Mathematical Ecology (Biomathematics 18), edited by Simon A. Levin, Thomas G. Hallam, Louis J. Gross, Springer, Berlin, 1989. 



\bibitem{KM} W. O. Kermack and A.  G. McKendrick, A contribution to the mathematical theory of epidemics, Proceedings of the Royal Society A, 115 (1927),  700--721.

\bibitem{M} Robert E. Mickens, An exactly solvable model for the spread of disease, College Mathematics Journal 43 (2012), 114--121. 

\end{thebibliography}
\end{document}